\documentclass[11pt]{amsart}
\usepackage{bm}
\usepackage{fullpage}
\usepackage{amssymb}
\usepackage{amsmath,amsfonts,amsthm}
\usepackage{hyperref}
\usepackage{color}
\usepackage{cite}
\usepackage{tikz}
\definecolor{uuuuuu}{rgb}{0.26666666666666666,0.26666666666666666,0.26666666666666666}

\newtheorem{theorem}{Theorem}[section]

\newtheorem{claim}[theorem]{Claim}
\newtheorem{proposition}[theorem]{Proposition}

\newtheorem{conjecture}{Conjecture}

\title{A PROOF OF A CONJECTURE ON RAMSEY NUMBERS $B(2,2,3)$ }
\author[1]{Yaser Rowshan$^1$}
 \author[1]{Mostafa  Gholami$^1$}

\address{$^1$Department of Mathematics, Institute for Advanced Studies in Basic Sciences (IASBS), Zanjan 66731-45137, Iran}

\email{y.rowshan@iasbs.ac.ir}
\email{gholami.m@iasbs.ac.ir}
 
\begin{document}
\maketitle

\begin{abstract}
 The bipartite Ramsey number $B(n_1,n_2,\ldots,n_t)$ is the least positive integer $b$ such that, any coloring of the edges of $K_{b,b}$ with
$t$ colors will result in a monochromatic copy of $K_{n_i,n_i}$ in the $i-$th color, for some $i$, $1\leq i\leq t$.  In this paper we obtain the exact values of bipartite Ramsey numbers $B(2,2,3)$. In particular, we prove  the conjecture of Radziszowski at al. aobut $B(2,2,3)$ which was introduced in 2015. In fact we prov that  $B(2,2,3)=17$.
\end{abstract}

\section{Introduction}
 The bipartite Ramsey number $B(n_1,n_2,\ldots,n_t)$ is the least positive integer $b$ such that, any coloring of the edges of $K_{b,b}$ with $t$ colors will result in a monochromatic copy of $K_{n_i,n_i}$ in the $i-$th color, for some $i$, $1\leq i\leq t$. The existence of such a positive integer is guaranteed by a result of Erd\H{o}s and Rado \cite{erdos1956partition}.\\
 The Zarankiewicz number $z(K_{m,n}, t)$ is defined as the maximum number of edges in any subgraph $G$ of the complete bipartite graph $K_{m,n}$, such that $G$ does not contain $K_{t,t}$ as a subgraph. Zarankiewicz numbers and related extremal graphs have been studied by numerous authors, including K\'{o}vari \cite{Kovari1954}, Reiman \cite{Reiman1958berEP} and  Goddard, Henning,
and Oellermann in \cite{goddard2000bipartite}.\\
The study of bipartite Ramsey numbers was initiated by Beineke and Schwenk in 1976 \cite{beinere1976bipartite}, and continued by others, in particular Exoo \cite{Exoo1991ABR}, Hattingh and Henning \cite{hattingh1998star}. 
The following exact values have been established.  $B(2,5)=17$ \cite{collins2016zarankiewicz},  $B(2,2,2,2)=19$ \cite{dybizbanski2015some}, $B(2,2,2)=11$ \cite{Exoo1991ABR}. In the smallest open case for $5$ colors, it is known that $26\leq B(2,2,2,2,2)\leq 28$  \cite{dybizbanski2015some}. One can refer to  \cite{rowshan2021size, raeisi2015star, hatala2021new, dybizbanski2015some, Kovari1954, collins2015bipartite, gholami2021bipartite} and it references for further studies.  Radziszowski et al.  in \cite{dybizbanski2015some} showed that $17\leq B(2,2,3)\leq18$ and in the same source has made the following conjecture:
\begin{conjecture}(\cite{dybizbanski2015some})
	$B(2,2,3)=17$.
\end{conjecture}
We intend to compute the exact values of the multicolor bipartite Ramsey numbers $B(2,2,3)$. Actually, we prove the following results: 

\bigskip
\begin{theorem}\label{th1} 
 $B(2,2,3)=17$.
\end{theorem}

In this paper, we only concerned with undirected, simple and finite graphs. We follow \cite{bondy1976graph} for terminology and notations not defined here. Let $G$ be a graph with vertex set $V(G)$ and edge set $E(G)$. The degree of a vertex $v\in V(G)$ is denoted by $\deg_G(v)$, or simply by $\deg(v)$. The neighborhood $N_G(v)$ of a vertex $v$ is the set of all vertices of $G$ adjacent to $v$ and satisfies $|N_G(v)|=\deg_G(v)$. The minimum and maximum degrees of vertices of $G$ are denoted by $\delta(G)$ and $\Delta(G)$, respectively.
Also, the complete bipartite graph with bipartition $(X,Y)$, where $|X|=m$ and $|Y|=n$ is denoted by $K_{m,n}$. We use $[X,Y]$ to denote the set of edges between a bipartition $(X,Y)$ of $G$. The degree of a vertex $v\in X\subseteq V(G)$ is denoted by $D_G(X)=(d_1,d_2,\ldots, d_{n})$ where $|X|=n$.
The complement of a graph $G$, denoted by $\overline{G}$, is the graph with the same vertices as $G$ and contains those edges which are not in $G$.
$G$ is $n$-colorable to $(G_1, G_2,\ldots, G_n)$ if there exist a $n$-edge decomposition of $G$, say $(H_1, H_2,\ldots, H_n)$ where $G_i\nsubseteq H_i$ for each $i=1,2, \ldots,n.$

\section{Proof of the main theorem}
To prove our main results, namely Theorem \ref{th1}, we need to establish some preliminary results. We start with the following proposition:


\begin{proposition}\label{pro1}{(\cite{collins2016zarankiewicz, collins2015bipartite})}The following results about Zarankiewicz number are true: 	
 	\begin{itemize}
 	\item[$\bullet$] $z(K_{17,17}, 2)=74$.
 	\item[$\bullet$] $z(K_{16,17}, 2)= 71$.
 	\item[$\bullet$] $z(K_{17,17}, 3)= 141$.
 	\item[$\bullet$] $z(K_{16,17}, 3)= 133$.
 	\item[$\bullet$] $z(K_{13,17}, 3)= 110$. 
 		\item[$\bullet$] $z(K_{12,17}, 3)= 103$. 
 	\item[$\bullet$] $z(K_{11,17}, 3)= 96$.
 
 \end{itemize}
\end{proposition}
\begin{proof}
	By using the bounds in Table $3$ and Table $4$ of \cite{dybizbanski2015some}  and  Table $C.3$ of \cite{collins2015bipartite} the  proposition holds.
\end{proof}
\begin{theorem}\label{t0}\cite{dybizbanski2015some}
 $17\leq B(2,2,3)\leq 18$.
\end{theorem}
\begin{proof}
The lower bound witness is found in Table $2$ of  \cite{dybizbanski2015some}. The upper bound is implied by using the bounds in Table $3$ and Table $4$ of \cite{dybizbanski2015some}. We know that $z(K_{18,18}, 2) = 81$, $z(K_{18,18}, 3) \leq 156$, and $2 \times 81 + 156 = 318 < 324=|E(K_{18,18})|$.
\end{proof}	

Suppose that $(G^r,G^b,G^g)$ be $3$-edges coloring  of $K_{17,17}$ where $K_{2,2} \nsubseteq G^r$, $K_{2,2} \nsubseteq G^b$ and $K_{3,3} \nsubseteq G^g$, in the following theorem, we specify some properties of the subgraph with color g. The properties regarding  $\Delta(G^g)$, $\delta (G^g)$, $E(G^g)$ and degree sequence of vertices $X$, $Y$ in induced graph with color $g$.

\begin{theorem}\label{t1}
	Assume that $(G^r,G^b,G^g)$ be $3$-edges coloring  of $K_{17,17}$ where $K_{2,2} \nsubseteq G^r$, $K_{2,2} \nsubseteq G^b$ and $K_{3,3} \nsubseteq G^g$. Hence  we have:
	\begin{itemize}
		\item[$~(a)$] $|E(G^g)|=141$,
		\item[(b)] $\Delta(G^g) =9$ and  $\delta(G^g) =8$,
		\item[(c)] $D_{G^g}(X)=D_{G^g}(Y)=(9,9,9,9,9,8,8,\ldots,8)$.
	\end{itemize}
	
\end{theorem}
\begin{proof}
Assume that  $X=\{x_1,x_2,\ldots,x_{17}\}$, $Y=\{y_1,y_2,\ldots,y_{17}\}$ be the partition set of $K=K_{17,17}$ and $(G^r,G^b,G^g)$ be $3$-edges coloring of $K$ where $K_{2,2} \nsubseteq G^r$, $K_{2,2} \nsubseteq G^b$ and $K_{3,3} \nsubseteq G^g$. Since $|E(K)|=289$, if $|E(G^g)|\leq 140$ then  $|E(\overline{G^g})|\geq 149$ that is, either $|E(G^r)|\geq 75$ or $|E(G^b)|\geq 75$. In any case by Proposition \ref{pro1}, either $K_{2,2} \subseteq G^r$ or $K_{2,2} \subseteq G^b$, a contradiction. Hence assume that $|E(G^g)|\geq 141$. If $|E(G^g)|\geq 142$ then by by Proposition \ref{pro1}, $K_{3,3}\subseteq G^g$, a contradiction again that is, we have $|E(G^g)|=141$ and the part $~(a)$ is true.\\
To prove the part $(b)$, by part $~(a)$ since $|E(G^g)|=141$, we can check that  $\Delta(G^g) \geq 9$. Assume that there exist a vertex of $V(K)$ say $x$ such that, $|N_{G^g}(x)|\geq 10$ that is, $\Delta(G^g) \geq 10$. Consider $x$, set $G_1^{g}=G^g\setminus \{x\}$, hence by part $~(a)$ we have $|E(G_1^{g})|\leq 141-10=131$ therefore, since $|E(K_{16,17})|=272$ we have $|E(\overline{G^{g}_1})|\geq  141$ that is, either $|E(G_1^{r})|\geq 71$ or $|E(G_1^{b})|\geq 71$. In any case by Proposition \ref{pro1} either $K_{2,2} \subseteq G_1^{r}\subseteq G^r$ or $K_{2,2} \subseteq G_1^{b}\subseteq G^b$, a contradiction. So, we have  $\Delta(G^g) =9$. To prove $\delta(G^g) =8$, assume that $M=\{x\in X, ~|N_{G^g}(x)|=9\}$ and  $N=\{x\in X, ~|N_{G^g}(x)|=8\}$, by part $~(a)$ one can say that $|M|\geq 5$, if $|M|= 6$ then, we have $\delta(G^g) \leq 7$ that is, there exist a vertex of $X$ say $x$ such that, $|N_{G^g}(x)|\leq 7$ therefore, $|N|\leq 10$. If $|N|=10$ then, we have $|E(G^g[M\cup N, Y])|=134$, so by  Proposition \ref{pro1} we have $K_{3,3}\subseteq G^g$, a contradiction. Now assume that $|N|\leq 9$ therefore,  $|E(G^g)|\leq (6\times 9)+ (9\times 8)+(2\times 7)=140$,  a contradiction again. For $|M|= 7$,  if $|N|\geq 6$, we have $|E(G^g[M\cup N', Y])|=111$ where $N'\subseteq N$ and $|N'|=6$ therefore, by  Proposition \ref{pro1} we have $K_{3,3}\subseteq G^g$, a contradiction, hence assume that $|N|\leq 5$ therefore, we have  $|E(G^g)|\leq (7\times 9)+ (5\times 8)+(5\times 7)=138$,  a contradiction again. For $|M|= 8$,  if $|N|\geq 5$ then, $|E(G^g[M\cup N', Y])|=112$ where $N'\subseteq N$ and $|N'|=5$ therefore, by  Proposition \ref{pro1} we have $K_{3,3}\subseteq G^g$, a contradiction, so assume that $|N|\leq 4$ that is,   $|E(G^g)|\leq (8\times 9)+ (4\times 8)+(5\times 7)=139$,  a contradiction again. For $|M|= 9$ if $|N|\geq 3$, then  $|E(G^g[M\cup N', Y])|=105$  where $N'\subseteq N$ and $|N'|=3$, so by  Proposition \ref{pro1} we have $K_{3,3}\subseteq G^g$, a contradiction, so $|N|\leq 2$ that is,  $|E(G^g)|\leq (9\times 9)+ (2\times 8)+(6\times 7)=139$  which is a contradiction again. For $|M|= 10$, if $|N|\geq 1$ then,  $|E(G^g[M\cup N', Y])|=98$  where $N'\subseteq N$ and $|N'|=1$ therefore, by  Proposition \ref{pro1}  $K_{3,3}\subseteq G^g$, a contradiction. Thus assume that $|N|=0$ so,  $|E(G^g)|\leq (10\times 9)+(7\times 7)=139$,  a contradiction again. Therefore,   $|M|=5$ and $|N|=12$ that is,  $\delta(G^g) =8$ and the part $~(b)$ is true.\\
Now by parts $~(a)$ and $(b)$ it is easy to say that  $D_{G^g}(X)=D_{G^g}(Y)=(9,9,9,9,9,8,8,\ldots,8)$ that is, the  part $~(c)$ is true and the proof is complete.  
\end{proof}
Suppose that $(G^r,G^b,G^g)$ be $3$-edge coloring  of $K_{17,17}$ where $K_{2,2} \nsubseteq G^r$, $K_{2,2} \nsubseteq G^b$ and $K_{3,3} \nsubseteq G^g$. In the following theorem we discuss about the maximmum number of common neighbors of $G^g(x)$ and $G^g(x')$ for $x, x' \in X$.
\begin{theorem}\label{t2}
	Assume that $(G^r,G^b,G^g)$ be $3$-edge coloring  of $K_{17,17}$ where $K_{2,2} \nsubseteq G^r$, $K_{2,2} \nsubseteq G^b$ and $K_{3,3} \nsubseteq G^g$. Let $|N_{G^g}(x)|=9$ and  $N_{G^g}(x)=Y_1$, the following results are true:
	\begin{itemize}
		\item[$~(a)$] For each $x\in X\setminus \{x_1\}$, we have $|N_{G^g}(x)\cap Y_1|\leq 5$,
		\item[(b)] Assume that $n=\sum \limits_{i=1}^{i=17}|N_{G^g}(x_i)\cap Y_1|$, then $72\leq n\leq 73$.
	\end{itemize}
\end{theorem}
\begin{proof} 
	Assume that  $X=\{x_1,x_2,\ldots,x_{17}\}$, $Y=\{y_1,y_2,\ldots,y_{17}\}$ be the partition set of $K=K_{17,17}$ and $(G^r,G^b,G^g)$ be $3$-edges coloring of $K$ where $K_{2,2} \nsubseteq G^r$, $K_{2,2} \nsubseteq G^b$ and $K_{3,3} \nsubseteq G^g$. Without loss of generality (W.l.g) assume that $x=x_1$ and $Y_1=\{y_1,\ldots,y_9\}$. To prove the part $~(a)$, by contrary  assume that there exist a vertex of $X\setminus \{x_1\}$ say $x$  such that, $|N_{G^g}(x)\cap Y_1|\geq 6$. W.l.g assume that $x=x_2$ and $Y_2=\{y_1,y_2,\ldots,y_6\}\subseteq N_{G^g}(x_2)$. Since $K_{3,3} \nsubseteq G^g$, for each $x\in X\setminus \{x_1,x_2\}$, we have $|N_{G^g}(x)\cap Y_2|\leq 2$ that is, $\sum \limits_{i=1}^{i=17}|N_{G^g}(x_i)\cap Y_2|\leq 6+6+(15\times 2)\leq 42$. Now, since $|E(G^g[X,Y_2])|\leq 42$, one can check that  there exist at least one vertex of $Y_2$ say $y$ such that,  $|N_{G^g}(y)|\leq 7$, a contradiction to part $(c)$ of Theorem \ref{t1}. Hence  $|N_{G^g}(x)\cap Y_1|\leq 5$ for each $x\in X\setminus \{x_1\}$ that is, the part $~(a)$ is true.\\
To prove the part $(b)$, if $n\leq 71$ then, by part $(c)$ of Theorem \ref{t1}, one can check that  there exist at least one vertex of $Y_1$ say $y$ such that,  $|N_{G^g}(y)|\leq 7$, a contradiction. Therefore,  $n\geq 72$. Assume that $n\geq 74$ and let $D_{G^g}(Y_1)=(d_1,d_2,\ldots,d_9)$. Since $\sum \limits_{i=1}^{i=17}|N_{G^g}(x_i)\cap Y_1|\geq 74$, there exist at least two vertices of $Y_1$ say $y', y''$ such that,  $|N_{G^g}(y')|=|N_{G^g}(y'')|=9$. Since  $n\geq 74$ and $|X\setminus \{x_1\}|=16$, there exist at lest one vertex of $X\setminus \{x_1\}$ say $x'$ such that, $|N_{G^g}(x')\cap Y_1|=5$.  W.l.g assume that  $x'=x_2$ and $N_{G^g}(x_2)\cap Y_1=Y_2=\{y_1,\ldots,y_5\}$. Now we note the following claims:
\begin{claim}\label{c1}
For each $x\in X\setminus \{x_1,x_2\}$ we have $|N_{G^g}(x)\cap Y_2|= 2$ and  $D_{G^g}(Y_2)=(8,8,8,8,8)$.
\end{claim}
Proof of the Claim. Since $K_{3,3} \nsubseteq G^g$ for each $x\in X\setminus \{x_1,x_2\}$ thus, $|N_{G^g}(x)\cap Y_2|\leq 2$ that is, $\sum \limits_{i=1}^{i=17}|N_{G^g}(x_i)\cap Y_2|\leq 5+5+(15\times 2)\leq 40$. Now, since $|E(G^g[X,Y_2])|\leq 40$ and $|Y_2|=5$, if there exist a vertex of $X_1$ say $x'$ such that, $|N_{G^g}(x)\cap Y_2|\leq 1$ then, $|E(G^g[X,Y_2])|\leq 39$ therefore, one can check that there exist at least one vertex of $Y_2$ say $y$ such that,  $|N_{G^g}(y)|\leq 7$, a contradiction to part $(c)$ of Theorem \ref{t1}. So  $|N_{G^g}(x)\cap Y_2|= 2$ and $\sum \limits_{y\in Y_2}|N_{G^g}(y)|=40$ therefor, by part $(c)$ of Theorem \ref{t1}  $D_{G^g}(Y_2)=(8,8,8,8,8)$ and the proof of the claim is compelete.
\begin{claim}\label{c2}
 $D_{G^g}(X_1)=(5,4,4,\ldots,4)$ where $X_1=X\setminus \{x_1\}$ in other word, $|N_{G^g}(x_i)\cap Y_1|=4$ for each $i\in\{3,4,\ldots,17\}$.
\end{claim}
Proof of the Claim. By contrary, assume that there exist a vertex of $X\setminus \{x_1,x_2\}$ say $x$ such that, $|N_{G^g}(x)\cap Y_1|=5$. W.l.g assume that $x=x_3$ and $N_{G^g}(x_3)\cap Y_1=Y_3$, now by Claim \ref{c1} we have $|N_{G^g}(x_3)\cap Y_2|=2$. W.l.g assume that $Y_3=\{y_1,y_2,y_6,y_7,y_8\}$ therefore, by Claim \ref{c1}, $D_{G^g}(Y_3)=(8,8,8,8,8)$ that is,  $|N_{G^g}(y)|=8$ for each $y\in Y_1\setminus \{y_9\}$. Since  $\Delta(G^g) =9$, we can check that   $n=\sum \limits_{i=1}^{i=17}|N_{G^g}(x_i)\cap Y_1|=\sum \limits_{i=1}^{i=9}|N_{G^g}(y_i)|\leq (8\times 8)+9=73$, a contradiction. So  $D_{G^g}(X_1)=(5,4,4,\ldots,4)$ and the proof of the claim is compelete.\\
Assume that $N_{G^g}(x_2)\cap Y_1=Y_2=\{y_1.\ldots,y_5\}$,  by Claim \ref{c1}  $D_{G^g}(Y_2)=(8,8,8,8,8)$. Since there exist at lest two vertices of $Y_1$ say $y',y''$ such that, $|N_{G^g}(y')|=|N_{G^g}(y'')|=9$ thus, $y',y''\in \{y_6,y_7,y_8,y_9\}$. W.l.g we may assume that $y'=y_6$ and $N_{G^g}(y_6)=X_2=\{x_1,x_3,\ldots, x_{10}\}$. By Claim \ref{c2} $|N_{G^g}(x)\cap Y_1|=4$ and  $|N_{G^g}(x)\cap Y_2|=2$ for each $x\in X_2\setminus \{x_1\}$ that is,  $|N_{G^g}(x)\cap \{y_7,y_8,y_9\}|=1$  for each $x\in X_2\setminus \{x_1\}$. Since $|X_2\setminus \{x_1\}|=8$ and $|N_{G^g}(x)\cap \{y_7,y_8,y_9\}|=1$, by the pigeon-hole principle we can check that there exist a vertex of $\{y_7,y_8,y_9\}$ say $y$ such that, $|N_{G^g}(y)\cap  X_2\setminus \{x_1\}|\geq 3$. W.l.g we may assume that $y=y_7$ and  $\{x_3,x_4,x_5\}\subseteq N_{G^g}(y_7)\cap  X_2\setminus \{x_1\}$. Since $|Y_2|=5$ and $|N_{G^g}(x_i)\cap Y_2|=2$ for $i=3,4,5$, we can say that there exist $i,i'\in \{3,4,5\}$ such that, $|N_{G^g}(x_i)\cap N_{G^g}(x_{i'})\cap  Y_2|\neq 0$, w.l.g assume that $i=3, i'=4$ and $y_1\in N_{G^g}(x_3)\cap N_{G^g}(x_{4})\cap  Y_2$. Therefore,  $K_{3,3}\subseteq G^g[\{x_1,x_3,x_4\},\{y_1,y_6,y_7\}]$, a contradiction. So,  $n\leq 73$ and the proof of the theorem is complete. 
	
\end{proof}
In the following theorem we prove that in any $3$-edge coloring  of $K_{17,17}$ say  $(G^r,G^b,G^g)$ where $K_{2,2} \nsubseteq G^r$, $K_{2,2} \nsubseteq G^b$, if there exist a vertex of $V(K)$ say $x$ such that, $|N_{G^g}(x)|=9$ and $\sum \limits_{x_i\in X\setminus\{x\}}|N_{G^g}(x_i)\cap N_{G^g}(x) |=64$ then, $K_{3,3} \subseteq G^g$.
\begin{theorem}\label{t3} 
	Assume that $(G^r,G^b,G^g)$ be $3$-edge coloring  of $K=K_{17,17}$ such that, $K_{2,2} \nsubseteq G^r$, $K_{2,2} \nsubseteq G^b$. Assume that there exist a vertex of $V(K)$ say $x$ such that, $|N_{G^g}(x)|=9$.  If  $~~\sum \limits_{i=1}^{i=17}|N_{G^g}(x_i)\cap Y_1|=73$ where $Y_1=N_{G^g}(x)$ then, $K_{3,3} \subseteq G^g$.
\end{theorem}
\begin{proof}
	By contrary, assume that $K_{3,3} \nsubseteq G^g$.  Therefore, by Theorem \ref{t1} and  Theorem \ref{t2} we have the following results:
\begin{itemize}
		\item[$(~a)$] $|E(G^g)|=141$, 
		\item[(b)] $\Delta(G^g) =9$ and  $\delta(G^g) =8$,
		\item[(c)] $D_{G^g}(X)=D_{G^g}(Y)=(9,9,9,9,9,8,8,\ldots,8)$,
		\item[(d)]For each $x'\in X\setminus \{x\}$ we have $|N_{G^g}(x)\cap N_{G^g}(x') |\leq 5$,	
		\item[(e)] If $A=\{x\in X,~|N_{G^g}(x)|=9\}$, then $|A|=5$ and  $72\leq\sum \limits_{y\in N_{G^g}(x)}|N_{G^g}(y)|\leq  73$, for each $x\in A$.
\end{itemize}
Assume that  $X=\{x_1,x_2,\ldots,x_{17}\}$, $Y=\{y_1,y_2,\ldots,y_{17}\}$ be the partition set of $K=K_{17,17}$ and $(G^r,G^b,G^g)$ be $3$-edges coloring of $K$ where $K_{2,2} \nsubseteq G^r$, $K_{2,2} \nsubseteq G^b$ and $K_{3,3} \nsubseteq G^g$. W.l.g assume that $x=x_1$, $Y_1=\{y_1,y_2,\ldots,y_9\}$ and $n=\sum \limits_{i=1}^{i=17}|N_{G^g}(x_i)\cap Y_1|=73$. Since $n=73$  by $(c)$ we can say that $D_{G^g}(Y_1)=(d_1,d_2,\ldots,d_9)=(9,8,8,\ldots,8)$ that is, there exist a vertex of $Y_1$ say $y$ such that,  $|N_{G^g}(y)|=9$. By (d) we have $|N_{G^g}(x_1)\cap N_{G^g}(x) |\leq 5$ for each $x \in X\{\setminus{x_1}\}$. Set $C=\{x\in X, ~|N_{G^g}(x)\cap N_{G^g}(x_1)|=5\}$. Now  by argument similar to proof of the Claim \ref{c1} we have the following claim:
	
\begin{claim}\label{cl3}
Assume that $x\in C$ and $N_{G^g}(x)\cap Y_1=Y'$ then, for each $x'\in X\setminus \{x_1,x\}$ we have $|N_{G^g}(x')\cap Y'|= 2$ and  $D_{G^g}(Y')=(8,8,8,8,8)$.
\end{claim}
Now we have the following claim about $|C|$:
		
\begin{claim}\label{cl4}
  $|C|\leq 2$.
\end{claim}
Proof of the Claim. By contrary, assume that $|C|\geq 3$. W.l.g assume that $\{x_2,x_3,x_4\}\subseteq C$ and $N_{G^g}(x_2)\cap Y_1=Y_2=\{y_1,\ldots, y_5\}$. By Claim \ref{cl3} $|N_{G^g}(x)\cap Y_2|= 2$ for each $x\in X\setminus \{x_1,x_2\}$.  W.l.g assume that   $N_{G^g}(x_3)\cap Y_1=Y_3=\{y_1,y_2,y_6, y_7, y_8\}$. Since $x_4 \in C$ and  $|N_{G^g}(x_4)\cap Y_i|= 2$ for $i=2,3$,  we have $y_9\in N_{G^g}(x_4)\cap Y_1$. Hence, for each $y\in Y_1$ there is at lest one $i\in \{2,3,4\}$ such that $y\in N_{G^g}(x_i)$ therefore, by Claim \ref{cl3} $D_{G^g}(Y_1)=(8,8,8,8,8,8,8,8,8)$, which is in contrast to   $\sum \limits_{i=1}^{i=17}|N_{G^g}(x_i)\cap Y_1|= \sum \limits_{i=1}^{i=9}|N_{G^g}(y_i)|=73$ so, $|C|\leq 2$.
	
\bigskip
Now by considering $|C|$ there are three cases as follow:
	
\bigskip
{\bf Case 1:} $|C|=0$. Since $n=73$, $|Y_1|=9$  and $|C|=0$ so, $D_{G^g}(X\setminus\{x_1\})=(4,4,\ldots,4,4)$,  $D_{G^g}(Y_1)=(9,8,8,8,8,,8,8,8,8)$,  $\sum \limits_{i=1}^{i=17}|N_{G^g}(x_i)\cap Y'|= 68$ and $D_{G^g}(Y')=(9,9,9,9,8,8,8,8)$ where $Y'=Y\setminus Y_1$.  Set $B=\{y\in Y',~~|N_{G^g}(y)|=9\}$ so,  $|B|=4$. \\
Now we are ready to prove the following claim:
\begin{claim}\label{cl5}
	There exist a vertex of $A\setminus \{x_1\}$ say $x$ such that:
	\[\sum_{y\in N_{G^g}(x)}|N_{G^g}(y)|\geq 74\] 
	Where $A=\{x\in X,  |N_{G^g}(x)|=9\}$.
\end{claim}
Proof of the Claim . Since $D_{G^g}(X_1)=(4,4,\ldots,4,4)$ and 	$D_{G^g}(Y_1)=(9,8,8,8,8,,8,8,8,8)$ for each $x\in A\setminus \{x_1\}$, one can say that: 
\[\sum_{y\in N_{G^g}(x)\cap Y_1}|N_{G^g}(y)|\geq 32\]
As $|Y'|=8$, $|B|=4$ and $|N_{G^g}(x_i)\cap Y'|= 5$ for each $x\in A\setminus \{x_1\}$, there exist at least one vertex of  $ A\setminus \{x_1\}$ say $x$ such that, $|N_{G^g}(x)\cap B|\geq 2$ otherwise, $K_{3,3}\subseteq G^g[A, Y'\setminus B]$, a contradiction. Hence w.l.g assume that $x_2\in A$ where $|N_{G^g}(x_2)\cap B|\geq 2$, so we have:
\[\sum_{y\in N_{G^g}(x_2)\cap Y'}|N_{G^g}(y)|\geq 42\]
That is, we have:
\[\sum_{y\in N_{G^g}(x_2)}|N_{G^g}(y)|=\sum_{y\in N_{G^g}(x_2)\cap Y'}|N_{G^g}(y)|+\sum_{y\in N_{G^g}(x)\cap Y_1}|N_{G^g}(y)|\geq 42+32= 74\] 
Now by considering  $x_2$ and $N_{G^g}(x_2)$ and by (e) (or part (b) of Theorem \ref{t2})  $K_{3,3}\subseteq G^g$, a contradiction again.\\
 	
\bigskip
{\bf Case 2:} $|C|=1$. W.l.g assume that $C=\{x_2\}$,  $N_{G^g}(x_2)\cap Y_1=Y_2=\{y_1,\ldots, y_5\}$. By Claim \ref{cl3} $|N_{G^g}(x_2)\cap N_{G^g}(x)\cap Y_1|=2$ for each $x\in X\setminus \{x_1,x_2\}$ and $|N_{G^g}(y_i)|=8$ for each $i\in \{1,2,\ldots,5\}$. Since there exist a vertex of $Y_1$ say $y$ such that,  $|N_{G^g}(y)|=9$, w.l.g we may assume that  $y=y_6$ and $N_{G^g}(y_6)=\{x_1, x_3,x_4\ldots, x_{10}\}$. Since $n=73$ and $|C|=1$ thus, $D_{G^g}(X_1)=(5,4,4,\ldots,4,3)$ that is, there exist at least seven  vertices of $N_{G^g}(y_6)\setminus \{x_1\}$ say $X_3=\{ x_3,x_4\ldots, x_{9}\}$ such that, $|N_{G^g}(x)\cap  Y_1|=4$ for each $x\in X_3$. Since $|X_3|= 7$, $|Y_2|=5$,  $|N_{G^g}(x)\cap  Y_1|=4$  and $|N_{G^g}(x_i)\cap Y_2|=2$ for each $x\in X_3$ thus, $|N_{G^g}(x)\cap  \{y_7,y_8,y_9\}|=1$ for each $x\in X_3$. Therefore, by the pigeon-hole principle we can say that there exist a vertex of $\{y_7,y_8,y_9\}$ say $y'$ such that, $|N_{G^g}(y')\cap  X_3|\geq 3$. W.l.g assume that $y'=y_7$ and $\{x_3,x_4,x_5\}\subseteq N_{G^g}(y_7)$. Therefore, since $|Y_2|=5$ one can check that there exist $i,i'\in \{3,4,5\}$ such that, $|N_{G^g}(x_i)\cap N_{G^g}(x_{i'})\cap  Y_2|\neq 0$.  W.l.g  assume that $i=3, i'=4$ and $y_1\in N_{G^g}(x_3)\cap N_{G^g}(x_{4})\cap  Y_2$. Therefore, $K_{3,3}\subseteq G^g[\{x_1,x_3,x_4\},\{y_1,y_6,y_7\}]$ which is  a contradiction.


\bigskip
{\bf Case 3:}$|C|=2$. W.l.g assume that $C=\{x_2,x_3\}$,  $N_{G^g}(x_2)\cap Y_1=Y_2=\{y_1,\ldots, y_5\}$. By Claim \ref{cl3} we have  $|N_{G^g}(x_2)\cap N_{G^g}(x_3)\cap Y_1|=2$. So, w.l.g we may assume that $N_{G^g}(x_3)\cap Y_1=Y_3=\{y_1,y_2,y_6,y_7, y_8\}$. Now by Claim \ref{cl3} we have $|N_{G^g}(y_i)|=8$ for each $i\in \{1,2,\ldots,8\}$. Since there is a vertex of $Y_1$ say $y$ such that, $|N_{G^g}(y)|=9$, we have $y=y_9$. W.l.g we may assume that $N_{G^g}(y_9)=X_2=\{x_1, x_4,x_5\ldots, x_{11}\}$. Since $n=73$ and $|C|=2$ thus, $D_{G^g}(X_1)=(5,5,4,4,\ldots,4,3,3)$ that is, there exist  two vertices of $X$ say $x,x'$ such that, $|N_{G^g}(x)\cap Y_1|=3$.  If $|N_{G^g}(y_9)\cap \{x,x'\}|\leq 1 $ then, there exist at least seven  vertices of $N_{G^g}(y_9)\setminus \{x_1\}$ such that, $|N_{G^g}(x)\cap  Y_1|=4$, in this case the proof is same as Case 1. Hence, assume that  $x,x' \in N_{G^g}(y_9)$. Since  $|N_{G^g}(x)\cap Y_2|=|N_{G^g}(x')\cap Y_2|=2$,  one can check that $|N_{G^g}(x)\cap \{y_6,y_7,y_8\}|=|N_{G^g}(x')\cap \{y_6,y_7,y_8\}|=0$.  Assume that $X_i=N_{G^g}(y_i)$ for $i=6,7,8$.  Since $|X_i|=8$ and $x,x'\notin X_i$ then,  for each  $x\in X_i\setminus \{x_1\}$ we have $|N_{G^g}(x)\cap  Y_1|=4$. Therefore,  by considering $X_i\setminus \{x_1\}$ and $y_i$ for each $i\in \{6,7,8\}$, the proof is same as  Case 1 and  $K_{3,3}\subseteq G^g$, a contradiction again.\\
 Therefore, by  Cases 1, 2, 3 the assumption dose not hold that is,   $K_{3,3}\subseteq G^g$ and the proof is complete.
\end{proof}	
In the following theorem we prove that in any $3$-edge coloring  of $K_{17,17}$ say  $(G^r,G^b,G^g)$ where $K_{2,2} \nsubseteq G^r$, $K_{2,2} \nsubseteq G^b$, if there exist a vertex of $V(K)$ say $x$ such that, $|N_{G^g}(x)|=9$ and $\sum \limits_{x_i\in X\setminus\{x\}}|N_{G^g}(x_i)\cap N_{G^g}(x) |=63$ then, $K_{3,3} \subseteq G^g$.
\begin{theorem}\label{t4} 
Assume that $(G^r,G^b,G^g)$ be $3$-edge coloring  of $K=K_{17,17}$ where $K_{2,2} \nsubseteq G^r$, $K_{2,2} \nsubseteq G^b$. Assume that there exist a vertex of $V(K)$ say $x$ such that, $|N_{G^g}(x)|=9$.  If  $~~\sum \limits_{i=1}^{i=17}|N_{G^g}(x_i)\cap Y_1|=72$ where $Y_1=N_{G^g}(x)$ then, $K_{3,3} \subseteq G^g$.
\end{theorem}
\begin{proof}
By contrary, assume that $K_{3,3} \nsubseteq G^g$.  Therefore, by Theorems \ref{t1} and \ref{t2} the following results are true:
\begin{itemize}
		\item[$(~a)$] $|E(G^g)|=141$,
		\item[(b)] $\Delta(G^g) =9$ and  $\delta(G^g) =8$,
		\item[(c)] $D_{G^g}(X)=D_{G^g}(Y)=(9,9,9,9,9,8,8,\ldots,8)$,	
		\item[(d)] For each $x\in X\setminus \{x_1\}$, we have $|N_{G^g}(x)\cap Y_1|\leq 5$,
	    \item[(e)] If $A=\{x\in X,~|N_{G^g}(x)|=9\}$ then, $|A|=5$ and  $72\leq\sum \limits_{y\in N_{G^g}(x)}|N_{G^g}(y)|\leq  73$, for each $x\in A$.
\end{itemize}
Assume that  $X=\{x_1,x_2,\ldots,x_{17}\}$, $Y=\{y_1,y_2,\ldots,y_{17}\}$ be the partition set of $K=K_{17,17}$ and $(G^r,G^b,G^g)$ be $3$-edges coloring of $K$ where $K_{2,2} \nsubseteq G^r$, $K_{2,2} \nsubseteq G^b$ and $K_{3,3} \nsubseteq G^g$. W.l.g assume that $x=x_1$, $Y_1=\{y_1,y_2,\ldots,y_9\}$ and $n=\sum \limits_{i=1}^{i=17}|N_{G^g}(x_i)\cap Y_1|=72$. Since $n=73$  by $(c)$ we can say that $D_{G^g}(Y_1)=(d_1,d_2,\ldots,d_9)=(8,8,8,\ldots,8)$. Set $C=\{x\in X, ~|N_{G^g}(x)\cap N_{G^g}(x_1)|=5\}$.  Define $D$ and $E$  as follow:
\[D=\{x\in X\setminus \{x_1\}, ~such~ that~  |N_{G^g}(x)\cap  Y_1|= 5\}\]
\[E=\{x\in X\setminus \{x_1\}, ~such~ that~  |N_{G^g}(x)\cap  Y_1|= 3\}\]
Now we have the following claim:
	
\begin{claim}\label{c6}
$|D|\leq 3$ and $|E|\leq 4$. 
\end{claim}
Proof of the Claim. By contrary, assume that $|D|\geq 4$. W.l.g  suppose that $\{x_2, x_3, x_4, x_5\}\subseteq D$,  $N_{G^g}(x_2)\cap Y_1=Y_2=\{y_1,\ldots, y_5\}$. Now, by  Claim \ref{cl3} we have  $|N_{G^g}(x)\cap Y_2|= 2$ for each $x\in X\setminus \{x_1, x_2\}$. W.l.g we may assume that  $N_{G^g}(x_3)\cap Y_1=Y_3=\{y_1,y_2,y_6, y_7, y_8\}$. Consider $N_{G^g}(x_i)\cap Y_1(i=4,5)$. Since $|N_{G^g}(x_i)\cap Y_j|= 2~(i=4,5,j=2,3)$ and $x_i\in A$ we can check that  $|N_{G^g}(x_i)\cap \{y_3,y_4,y_5 \}|= 2$, $|N_{G^g}(x_i)\cap \{y_6,y_7,y_8 \}|= 2$  and $y_9\in N_{G^g}(x_i)$ for $i=4,5$ otherwise, if there exist a vertex of $\{x_4,x_5\}$ say $x$ such that,  $|N_{G^g}(x_i)\cap \{y_1,y_2\}|\neq 2$ then,  one can say that $K_{3,3}\subseteq G^g[\{x_1,x_i,x\},Y_1]$ for some $i\in \{1,2\}$, a contradiction. Therefore, since $|\{y_3,y_4,y_5 \}|=|\{y_6,y_7,y_8 \}|=3$ and $x_4, x_5 \in A$, by the pigeon-hole principle $|N_{G^g}(x_4)\cap N_{G^g}(x_5)\cap  \{y_3,y_4,y_5 \}|\geq 1$ and $|N_{G^g}(x_4)\cap N_{G^g}(x_5)\cap  \{y_6,y_7,y_8 \}|\geq 1$. W.l.g we may assume that $y_3, y_6 \in N_{G^g}(x_4)\cap N_{G^g}(x_5)$ therefore, since $y_9\in N_{G^g}(x_4)\cap N_{G^g}(x_5)$,  we have $K_{3,3}\subseteq G^g[\{x_1,x_4,x_5\},\{y_3,y_6,y_9\}]$, a contradiction. Therefore, $|D|\leq 3$.  Now, as $\sum \limits_{i=2}^{i=17}|N_{G^g}(x_i)\cap Y_1|=63$ and $|D|\leq 3$,  we can say that $|E|\leq 4$ and the proof of claim is complete.\\ 
Now, by considering $|D|$, there are three cases as follow:

	\bigskip
{\bf Case 1:} $|D|=0$. Since $n=72$ and $|D|=0$ thus? $D_{G^g}(X\setminus \{x_1\})=(4,4,\ldots,4,3)$  and $D_{G^g}(Y_1)=(8,8,8,8,8,8,8,8,8)$,  $\sum \limits_{i=1}^{i=17}|N_{G^g}(x_i)\cap Y'|= 69$ and $D_{G^g}(Y')=(9,9,9,9,9,8,8,8)$ where $Y'=Y\setminus Y_1$. Set $B=\{y\in Y',  |N_{G^g}(y)|=9\}$ Hence, $|B|=5$.\\ 
Now, we have the following claim:

\begin{claim}\label{cl5}
	There exist a vertex of $A\setminus \{x_1\}$ say $x$ such that:
	\[\sum_{y\in N_{G^g}(x)}|N_{G^g}(y)|\geq 75\]
	Where $A=\{x\in X,  |N_{G^g}(x)|=9\}$.
\end{claim} 
Proof of the Claim. Since $D_{G^g}(X_1)=(4,4,\ldots,4,3)$ and 	$D_{G^g}(Y_1)=(8,8,8,8,8,,8,8,8,8)$ so, for at least three vertices of  $A\setminus \{x_1\}$ one can say that: 
\[\sum_{y\in N_{G^g}(x)\cap Y_1}|N_{G^g}(y)|\geq 32\]
Therefore, since $|N_{G^g}(x_i)\cap Y'|= 5$ for each $x\in A\setminus \{x_1\}$ and $D_{G^g}(Y')=(9,9,9,9,9,8,8,8)$, there exist at least one vertex of  $  A\setminus \{x_1\}$ say $x$ such that, $|N_{G^g}(x)\cap B|\geq 3$ otherwise, $K_{3,3}\subseteq G^g[A, Y'\setminus B]$, a contradiction. Hence, w.l.g assume that $x_2\in A$ and $|N_{G^g}(x_2)\cap B|\geq 3$ therefore, we have:
\[\sum_{y\in N_{G^g}(x)\cap Y''}|N_{G^g}(y)|\geq 3\times9 + 2\times 8= 43\]
That is, we have:
\[\sum_{y\in N_{G^g}(x_2)}|N_{G^g}(y)|=\sum_{y\in N_{G^g}(x_2)\cap Y'}|N_{G^g}(y)|+\sum_{y\in N_{G^g}(x)\cap Y_1}|N_{G^g}(y)|\geq 43+32= 75\] 
Now by considering  $x_2$ and $N_{G^g}(x_2)$ and by $(e)$( or by part $(b)$ of Theorem \ref{t2}), $K_{3,3}\subseteq G^g$, a contradiction again. 	
		 
\bigskip
{\bf Case 2:} $|D|= 1$ (for the case that $(|D|= 2)$ the proof is same). W.l.g assume that $D=\{x_2\}$, $N_{G^g}(x_2)\cap Y_1=Y_2=\{y_1,\ldots, y_5\}$. Since $n=72$, $|D|= 1$ and $|N_{G^g}(x)\cap  Y_1|\leq 5$ we can say that $|E|=2$. As $|N_{G^g}(x)\cap Y_2|= 2$ for each $x\in X\setminus \{x_1, x_2\}$ and $|E|=2$, we can check that there exist a vertex of $\{y_6,y_7,y_8,y_9\}$ say $y$ such that, for each vertex of $N_{G^g}(y)\cap X\setminus \{x_1\}$ say $x$, $|N_{G^g}(x)\cap Y_1|= 4$. W.l.g we may assume that $y=y_6$, $N_{G^g}(y_6)\cap X\setminus \{x_1\}=\{x_3,x_4,\ldots, x_9\}$. Since $|N_{G^g}(y_6)\cap X\setminus \{x_1\}|=7$ and $|N_{G^g}(x)\cap Y_2|= 2$ for each $x\in N_{G^g}(y_6)\cap X\setminus \{x_1\}$ so, $|N_{G^g}(x)\cap \{y_7,y_8,y_9\}|= 1$. Therefore, by the pigeon-hole principle we can check that there exist a vertex of  $\{y_7,y_8,y_9\}$ say $y'$ such that, $|N_{G^g}(y_6)\cap N_{G^g}(y')\cap X\setminus \{x_1\}|\geq 3$. W.l.g assume that $y'=y_7$ and $\{x_3,x_4,x_5\}\subseteq N_{G^g}(y_6)\cap N_{G^g}(y_7)\cap X\setminus \{x_1\}$. Therefore, since $|Y_2|=5$ and $|N_{G^g}(x)\cap  Y_2|= 2$, there exist at least two vertices of $\{x_3,x_4,x_5\}$ say $x',x''$ such that, $|N_{G^g}(x')\cap N_{G^g}(x'')\cap Y_2|\neq 0$. W.l.g assume that $x'=x_3,x''=x_4$ and $y_1\in N_{G^g}(x_3)\cap N_{G^g}(x_4)$ therefore, we have $K_{3,3}\subseteq G^g[\{x_1,x_3,x_4\},\{y_1,y_6,y_7\}]$, a contradiction.

\bigskip
{\bf Case 3:} $|D|=3$. W.l.g assume that $D=\{x_2,x_3,x_4\}$,  $N_{G^g}(x_2)\cap Y_1=Y_2=\{y_1,\ldots, y_5\}$. By Claim \ref{cl3} we have $|N_{G^g}(x_2)\cap N_{G^g}(x_3)\cap Y_1|=2$. Wl.g we may assume that $N_{G^g}(x_3)\cap Y_1=Y_3=\{y_1, y_2, y_6,y_7, y_8\}$. Since $x_4\in D$ and  $|N_{G^g}(x_4)\cap Y_i|=2$ for $i=2,3$ thus, $y_9\in N_{G^g}(x_4)$. If $|N_{G^g}(x_4)\cap \{y_1,y_2\}|\neq 0$,  as $|N_{G^g}(x_2)\cap N_{G^g}(x_4)\cap Y_1|=2$ and $x_4\in D$, one can check that $|N_{G^g}(x_4)\cap \{y_6,y_7,y_8\}|=2$ that is, $K_{3,3}\subseteq G^g[\{x_1,x_3,x_4\},Y_1]$, a contradiction. Hence, $|N_{G^g}(x_4)\cap \{y_1,y_2\}|= 0$ therefore, $|N_{G^g}(x_4)\cap \{y_3,y_4,y_5\}|=2$  and $|N_{G^g}(x_4)\cap \{y_6,y_7,y_8\}|=2$. W.l.g we may assume that $N_{G^g}(x_4)\cap Y_1=Y_4=\{y_3, y_4, y_6,y_7, y_9\}$. Since $|D|=3$ we can say that $|E|=4$. W.l.g assume that $E=\{x_5,x_6,x_7,x_8\}$. Now, we have the following claim:
	
\begin{claim}\label{c8}
$|N_{G^g}(y_9)\cap E|=0$. 
\end{claim}
Proof of the Claim. By contrary, assume that $|N_{G^g}(y_9)\cap E|\neq 0$.  Assume that $x_5\in N_{G^g}(y_9)\cap E$ that is, $x_5y_9\in E(G^g)$.  Since $x_5\in E$ and  $\{x_2,x_3, x_4\}= D$,  by Claim \ref{cl3} we have  $|N_{G^g}(x_5)\cap N_{G^g}(x_i)|=|N_{G^g}(x_5)\cap Y_i|=2$ for $i=2,3,4$. Consider  $N_{G^g}(x_5)\cap Y_2$, assume that $N_{G^g}(x_5)\cap Y_2= \{y',y''\}$, if $\{y',y''\}=\{y_1,y_2\}$ then, $|N_{G^g}(x_5)\cap Y_4|=1$, a contradiction. Therefore, we may assume that $|\{y',y''\}\cap \{y_1,y_2\}|\leq 1$. If $|\{y',y''\}\cap \{y_1,y_2\}|=0$  then,  $|N_{G^g}(x_5)\cap Y_3|=0$ and if $|\{y',y''\}\cap \{y_1,y_2\}|=1$ then,  $|N_{G^g}(x_5)\cap Y_3|\leq 1$, in any case there exist a vertex of $D$ say $x'$ such that, $|N_{G^g}(x_5)\cap N_{G^g}(x')|=1$,  a contradiction. So, the assumption dose not hold and the claim is true.\\
Therefore, by Claim \ref{c8}, since $|N_{G^g}(y_9)\cap D|=0$, we can say that for any vertex  of $N_{G^g}(y_9)\cap X\setminus \{x_1\}$ say $x$, $|N_{G^g}(x)\cap Y_1|\geq 4$ therefore, by considering  $Y_2$ and $y_9$, as   $|N_{G^g}(y_9)\cap X\setminus \{x_1\}|=7$ and  $|N_{G^g}(x)\cap Y_1|\geq 4$ for each $x\in N_{G^g}(y_9)\cap X\setminus \{x_1\}$  the proof is same as Case 1, a contradiction.\\
 Therefore, by  Cases 1, 2, 3 the assumption dose not hold that is,   $K_{3,3}\subseteq G^g$ and the proof is complete.	
\end{proof}
 
 \bigskip
 Now, combining Theorems \ref{t1}, \ref{t2} \ref{t3} and \ref{t4}, yield the  proof of the Theorem \ref{th1}.
\bibliographystyle{spmpsci} 
\bibliography{BI}

\begin{thebibliography}{10}
\providecommand{\url}[1]{{#1}}
\providecommand{\urlprefix}{URL }
\expandafter\ifx\csname urlstyle\endcsname\relax
  \providecommand{\doi}[1]{DOI~\discretionary{}{}{}#1}\else
  \providecommand{\doi}{DOI~\discretionary{}{}{}\begingroup
  \urlstyle{rm}\Url}\fi

\bibitem{beinere1976bipartite}
BEINERE, L., LW, B., AJ, S.: On a bipartite form of the ramsey problem.  (1976)

\bibitem{bondy1976graph}
Bondy, J.A., Murty, U.S.R., et~al.: Graph theory with applications, vol. 290.
\newblock Macmillan London (1976)

\bibitem{collins2015bipartite}
Collins, A.F.: Bipartite Ramsey Numbers and Zarankiewicz Numbers.
\newblock Rochester Institute of Technology (2015)

\bibitem{collins2016zarankiewicz}
Collins, A.F., Riasanovsky, A.W., Wallace, J.C.: Zarankiewicz numbers and
  bipartite ramsey numbers.
\newblock Journal of Algorithms and Computation \textbf{47}, 63--78 (2016)

\bibitem{dybizbanski2015some}
Dybizba{\'n}ski, J., Dzido, T., Radziszowski, S.: On some zarankiewicz numbers
  and bipartite ramsey numbers for quadrilateral.
\newblock ARS COMBINATORIA \textbf{119}, 275--287 (2015)

\bibitem{erdos1956partition}
Erd{\"o}s, P., Rado, R.: A partition calculus in set theory.
\newblock Bulletin of the American Mathematical Society \textbf{62}(5),
  427--489 (1956)

\bibitem{Exoo1991ABR}
Exoo, G.: A bipartite ramsey number.
\newblock Graphs and Combinatorics \textbf{7}, 395--396 (1991)

\bibitem{gholami2021bipartite}
Gholami, M., Rowshan, Y.: The bipartite ramsey numbers $br(c_8, c_{2n})$ (2021)

\bibitem{goddard2000bipartite}
Goddard, W., Henning, M.A., Oellermann, O.R.: Bipartite ramsey numbers and
  zarankiewicz numbers.
\newblock Discrete Mathematics \textbf{219}(1-3), 85--95 (2000)

\bibitem{hatala2021new}
Hatala, I., H{\'e}ger, T., Mattheus, S.: New values for the bipartite ramsey
  number of the four-cycle versus stars.
\newblock Discrete Mathematics \textbf{344}(5), 112320 (2021)

\bibitem{hattingh1998star}
Hattingh, J.H., Henning, M.A.: Star-path bipartite ramsey numbers.
\newblock Discrete Mathematics \textbf{185}(1-3), 255--258 (1998)

\bibitem{Kovari1954}
Kóvari T., S.V.T.P.: On a problem of k. zarankiewicz.
\newblock Colloquium Mathematicae \textbf{3}(1), 50--57 (1954).
\newblock \urlprefix\url{http://eudml.org/doc/210011}

\bibitem{raeisi2015star}
Raeisi, G.: Star-path and star-stripe bipartite ramsey numbers in
  multicoloring.
\newblock Transactions on Combinatorics \textbf{4}(3), 37--42 (2015)

\bibitem{Reiman1958berEP}
Reiman, I.: {\"U}ber ein problem von k. zarankiewicz.
\newblock Acta Mathematica Academiae Scientiarum Hungarica \textbf{9}, 269--273
  (1958)

\bibitem{rowshan2021size}
Rowshan, Y., Gholami, M., Shateyi, S.: The size, multipartite ramsey numbers
  for nk2 versus path--path and cycle.
\newblock Mathematics \textbf{9}(7), 764 (2021)

\end{thebibliography}
\end{document}